\documentstyle{article}
\begin{document}
\baselineskip+8pt
\small \begin{center}{\textit{ In the name of
 Allah, the Beneficent, the Merciful}}\end{center}

\begin{center}{\bf  MATRIX REPRESENTATIONS FOR SYMMETRIC AND ANTISYMMETRIC MULTI-LINEAR  MAPS}
\end{center}
\begin{center}{\bf  Ural Bekbaev  \\Turin Polytechnic University in Tashkent,\\
 INSPEM, Universiti Putra Malaysia.\\
e-mail: bekbaev@science.upm.edu.my}
 \end{center}

\begin{abstract}{In this paper the main results in arXiv0901.3179v3, related to the matrix representation of
polynomial maps, are restated in traditional way of linear algebra assuming that variable vectors are
presented as column vectors. Some new results related to that subject are also included. Here one can
find the behavior of the matrices of polynomial maps with respect to the change of variables (coordinate system),
 matrix representations for symmetric and antisymmetric multi-linear maps. It is shown also that
 the offered representations are in good concordance with known operations over such multi-linear maps.  }\end{abstract}


  \vspace{0.5cm}

 Matrix representation for polynomial maps is offered in [1]
 and its application to find a radius of absolute convergence of power series in many variables is demonstrated in
 [2].
 In those papers vectors with given variable components
 are considered as row vectors which is not traditional in linear algebra.
 As far as now we are dealing with matrix representations of nonlinear maps to follow the tradition of linear algebra
 variable vectors should be considered as column vectors.
 Therefore in this paper we reformulate some main results of those papers in traditional way,
 including some new results, and offer matrix representations for symmetric and antisymmetric multi-linear maps.

 For a positive integer $n$ let $I_n$ stand for all row $n$-tuples with nonnegative integer entries with the
 following linear order: $\beta=
 (\beta_1,\beta_2,...,\beta_n)<\alpha=(\alpha_1,\alpha_2,...,\alpha_n)$ if and only if $\vert \beta\vert < \vert \alpha\vert$ or
$\vert \beta\vert = \vert \alpha\vert$ and $\beta_1> \alpha_1$ or
$\vert \beta\vert = \vert \alpha\vert$, $\beta_1= \alpha_1$ and
$\beta_2> \alpha_2$ etcetera, where $\vert \alpha\vert$ stands for
$\alpha_1+\alpha_2+...+\alpha_n $.

It is clear that for $\alpha, \beta, \gamma \in I_n$ one has
$\alpha < \beta $ if and only if $\alpha + \gamma < \beta
+\gamma$. We write $\beta \ll \alpha$ if $\beta_i \leq \alpha_i$
for all $i=1,2,...,n$, $\left(\begin{array}{c}
  \alpha \\
  \beta \\
\end{array}\right)$ stands for $\frac{\alpha!}{\beta!
(\alpha -\beta)!}$, $\alpha !=\alpha_1!\alpha_2!...\alpha_n!$.

 In future $n$, $n'$ and $n''$ are assumed to be any fixed nonnegative integers (In the case of $n=0$ it is assumed that $I_n=  \{0\}$).

 For any commutative, associative ring $R$ containing the field of rational numbers and nonnegative
integer numbers $p,p'$ let  $M_{n,n'}(p,p';R)=M(p,p';R)$ stand for
all $"p\times p'"$ size matrices $A=(A_{\alpha,\alpha'})_{\vert
\alpha \vert=p, \vert \alpha'\vert=p'}$ ($\alpha$ presents
row,$\alpha'$ presents column and $\alpha\in I_n,\alpha'\in
I_{n'}$) with entries from $R$. The ordinary size of a such matrix
is $\left(\begin{array}{c}
  p+n-1 \\
  n-1 \\
\end{array}\right)\times\left(\begin{array}{c}
  p'+n'-1 \\
  n'-1 \\
\end{array}\right)$. Over such kind matrices in addition to the ordinary sum and
product of matrices we consider the following "product" as well:

{\bf Definition 1.} If $A\in M(p,p';R)$ and $B\in M(q,q';R)$ then
$A\bigodot B=C\in M(p+q, p'+q';R)$  such that for any
$\vert\alpha\vert=p+q$, $\vert\alpha'\vert=p'+q'$, where
$\alpha\in I_n,\alpha'\in I_{n'}$,
$$C_{\alpha,\alpha'}=\sum_{\beta,\beta'}\left(\begin{array}{c}
  \alpha' \\
  \beta' \\
\end{array}\right)
    A_{\beta,\beta'}B_{\alpha-\beta,\alpha'-\beta'}
$$, where the sum is taken over all $\beta\in I_n,\beta'\in I_{n'}$, for which $\vert
\beta\vert=p$, $\vert \beta'\vert=p'$, $\beta\ll \alpha$ and
$\beta'\ll \alpha'$.

{\bf Example 1.} If $n=n'=2$, $p=p'=q=q'=1$ in ordinary notations
for matrices $A=\left(
\begin{array}{cc}
  a_{11} & a_{12} \\
  a_{21} & a_{22} \\
\end{array}
\right)$, $B=\left(
\begin{array}{cc}
  b_{11} & b_{12} \\
  b_{21} & b_{22} \\
\end{array}
\right)$ the product $A\bigodot B$ can be given as
$$A\bigodot B=\left(
\begin{array}{ccc}
  2a_{11}b_{11} & a_{11}b_{12}+a_{12}b_{11} & 2a_{12}b_{12} \\
  2(a_{11}b_{21}+a_{21}b_{11}) & a_{11}b_{22}+a_{22}b_{11}+a_{12}b_{21}+a_{21}b_{12} & 2(a_{12}b_{22}+a_{22}b_{12}) \\
  2a_{21}b_{21} & a_{21}b_{22}+a_{22}b_{21} & 2a_{22}b_{22} \\
\end{array}
\right)$$

Let us agree that $h$ ($H$, $v$, $V$) stands for any element of
$M(0,1;R)$ (respect. $M(0,p;R)$, $M(1,0;R)$, $M(p,0;R)$ , where
$p$ may be any nonnegative integer). For the sake of convenience
it will be assumed that $A_{\alpha,\alpha'}=0$ ($\alpha !=\infty$)
whenever $\alpha \notin I_n$ or $\alpha' \notin I_{n'}$ (respect.
$\alpha \notin I_n$).

{\bf Proposition 1.} For the above defined product the following
are true.

1. $A\bigodot B=B\bigodot A$.

2. $(A+B)\bigodot C=A\bigodot C+ B\bigodot C$.

3. $(A\bigodot B)\bigodot C=A\bigodot (B\bigodot C)$

4. $ (\lambda A)\bigodot B=\lambda (A\bigodot B)$ for any
$\lambda\in R$

5. If $A$ and $B$ are square upper triangular matrices then $A\bigodot B $
is also an upper triangular matrix.

6. If $R$ is an integral domain then $A\bigodot B=0$ if and only
    if $A=0$ or $B=0$.

7. $(A\bigodot V)B=(AB)\bigodot V$

8. $A(B\bigodot H)=(AB)\bigodot H$

In future $A^{(m)}$ means the $m$-th power of matrix $A$ with
respect to the new product.

{\bf Proposition 2.} If $h=(h_1,h_2,...,h_{n'})\in M(0,1;R)$,
$v=(v_1,v_2,...,v_n)\in M(1,0;R)$, then
$$(h^{(m)})_{0,\alpha'}=m!h^{\alpha'}, \hspace{1cm} (v^{(m)})_{\alpha,0}=\left(\begin{array}{c}
  m \\
  \alpha\\
\end{array}\right)v^{\alpha}$$, where $v^{\alpha}$ stands for
$v_1^{\alpha_1}v_2^{\alpha_2}...v_n^{\alpha_n}$

{\bf Proposition 3.} For any nonnegative integers $p$, $q$, $p'$,
$q'$ and matrices $A\in M_{n,n'}(p,p';R)$, $B\in
M_{n,n'}(q,q';R)$, $h=(h_1,h_2,...,h_n)\in M_{n,n}(0,1;R)$,
$v=(v_1,v_2,...,v_{n'})\in M_{n',n'}(1,0;R)$ the following
equalities
$$(\frac{h^{(p)}}{p!}A)\bigodot
(\frac{h^{(q)}}{q!}B)=\frac{h^{(p+q)}}{(p+q)!}(A\bigodot B),
\hspace{1cm} (A\frac{v^{(p')}}{p'!})\bigodot
(B\frac{v^{(q')}}{q'!})=(A\bigodot
B)\frac{v^{(p'+q')}}{(p'+q')!}$$
 are true.

 In future let us assume that $R$ is a field of characteristic
 zero. For every vector space considered it is supposed that some
 basis in it is fixed and elements of the vector space are written
 as column vectors.

 {\bf Corollary 1.} If $A_1,A_2,...,A_k$ are $"1\times 1"$ size
 square matrices and $v$ is a common eigenvector for all
 $A_1,A_2,...,A_k$ with eigenvalues
 $\lambda_1,\lambda_2,...,\lambda_k$, respectively, then
 $v^{(k)}$ is an eigenvector for $A_1\bigodot A_2\bigodot ...\bigodot
 A_k$ with the eigenvalue $k!\lambda_1\lambda_2...\lambda_k.$

{\bf Corollary 2.} The equality $Av^1\bigodot Av^2\bigodot
...\bigodot Av^m=A^{(m)}\frac{v^1\bigodot v^2\bigodot ... \bigodot
v^{m}}{m!}$ is true, where $A\in M(1,1;R)$.

{\bf Proof.} To prove Corollary 2 one can consider the equality
$$A_1(v^1+ v^2+ ... +v^m)\bigodot A_2(v^1+ v^2+ ... +v^m)\bigodot...\bigodot A_k(v^1+ v^2+ ... +v^m)
=$$ $$(A_1\bigodot A_2\bigodot ...\bigodot
 A_k)\frac{(v^1+ v^2+ ... +v^m)^{(k)}}{k!}$$
and compare the corresponding terms to get equality
$$\sum_{P}(\bigodot_{i\in P_1}A_iv^1)\bigodot(\bigodot_{i\in P_2}A_iv^2)\bigodot ... \bigodot (\bigodot_{i\in
P_m}A_iv^m)=$$ $$(A_1\bigodot A_2\bigodot ...\bigodot
 A_k)(\frac{(v^1)^{(\alpha_1)}}{\alpha_1!}\bigodot
\frac{(v^2)^{(\alpha_2)}}{\alpha_2!}\bigodot ... \bigodot
\frac{(v^m)^{(\alpha_m})}{\alpha_m!})$$ , where
$\alpha_1+\alpha_2+...+\alpha_m=k$ and the sum is taken over all
partitions $P=(P_1, P_2,...,P_m)$ of $\{1, 2, ..., k\}$ for
which $|P_1|=\alpha_1$, $|P_2|=\alpha_2$,..., $|P_m|=\alpha_m$.
The Corollary is a particular case of this equality, when $m=k$,
$\alpha_1= \alpha_2=...=\alpha_m=1$, $A=A_1=A_2=...=A_k$.

Here are some relations between properties of $A \in M(1,1;R)$ and
$\frac{A^{(k)}}{k!}$.

{\bf Theorem 1.} Let $A\in M(1,1;R)$ be a square matrix.

1. If $\{\lambda_1, \lambda_2,..., \lambda_n\}$ is its all eigenvalues, where
each eigenvalue occurs as many times as its multiplicity, then
$\{\lambda_1^{\alpha_1}\lambda_2^{\alpha_2}...\lambda_n^{\alpha_n}\}_{|\alpha|=k}$
represents all eigenvalues of $\frac{A^{(k)}}{k!}$ including
their multiplicities as eigenvalues of  $\frac{A^{(k)}}{k!}$.
Moreover if $Av_i=\lambda_iv_i$ for $i=1,2,...,n$ then
$$\frac{A^{(k)}}{k!}(v_1^{(\alpha_1)}\bigodot
v_2^{(\alpha_2)}\bigodot...\bigodot
v_n^{(\alpha_n)})=\lambda_1^{\alpha_1}\lambda_2^{\alpha_2}...\lambda_n^{\alpha_n}(v_1^{(\alpha_1)}\bigodot
v_2^{(\alpha_2)}\bigodot...\bigodot v_n^{(\alpha_n)})$$

2. If $rk(A)=l$ then $$rk(\frac{A^{(k)}}{k!})=\left(\begin{array}{c}
  k+l-1 \\
  l-1 \\
\end{array}\right)$$

3. $$\det \frac{A^{(k)}}{k!}=(\det A)^{\left(\begin{array}{c}
  k+n-1 \\
  n \\
\end{array}\right)}$$

{\bf Proof.} A proof of this Theorem can be derived from the equality
$$\frac{(T^{-1}AT)^{(k)}}{k!}=\frac{(T^{(k)})^{-1}}{k!}\frac{A^{(k)}}{k!}\frac{T^{(k)}}{k!}$$,
where $T\in M(1,1;R)$. Note that $\frac{(T^{(k)})^{-1}}{k!}$ is
inverse for $\frac{T^{(k)}}{k!}$.

Indeed by choose of $T$ one can make $T^{-1}AT$ in Jordan form. In
common case $\frac{(T^{-1}AT)^{(k)}}{k!}$ will not be in Jordan
form but it will be in upper triangular form with
$\lambda_1^{\alpha_1}\lambda_2^{\alpha_2}...\lambda_n^{\alpha_n}$
on the main diagonal, where $|\alpha|=k$. Therefore the first
statement of Theorem 1 is true.

The proofs of other statements of Theorem 1 can be done in a similar way.

 In future the expression $Exp(A)$, whenever it has meaning, stands
for
$$E+\frac{1}{1!}A+\frac{1}{2!}A^{(2)}+\frac{1}{3!}A^{(3)}+ . .
.=\sum_{i=0}^{\infty}\frac{1}{i!}A^{(i)}$$, $R[x]$ is the ring of
polynomials in variables $x_1,x_2,...,x_{n'}$  over $R$,
$x=(x_1,x_2,...,x_{n'})\in M_{n',n'}(1,0;R[x])$.

If $n\neq 0$ and
$$\varphi(x)=(\varphi_1(x),\varphi_2(x),...,\varphi_{n}(x))=$$
$$M_{\varphi}(1,0)x^{(0)}+M_{\varphi}(1,1)\frac{x^{(1)}}{1!}+
M_{\varphi}(1,2)\frac{x^{(2)}}{2!}+...\in M_{n',n}(0,1;R)$$ is a
polynomial map from $R^{n'}$ to $R^n$ then
 one can screen it in the form
$$\varphi(x)=M_{\varphi}Exp{(x)} $$, where $M_{\varphi}\in
Mat(R)$ with blocks $M_{\varphi}(p,p')$ such that
$M_{\varphi}(p,p')=0$ whenever $p\neq 1$ and only finite number
blocks of the form $M_{\varphi}(1,p')$ are not zero, $Mat(R)$
stands for the set of all block-matrices of the form
$A=(A(p,p'))_{p=\overline{0,\infty},p'=\overline{0,\infty}}$ with
$A(p,p')\in M(p,p';R)$. We call $M_{\varphi}$ the matrix of the
polynomial map $\varphi(x)$. (Of course, if $n'=0$ then
$$\varphi(x)=M_{\varphi}(0,0)x^{(0)}+M_{\varphi}(0,1)\frac{x^{(1)}}{1!}+
M_{\varphi}(0,2)\frac{x^{(2)}}{2!}+...\in M_{0,n}(0,0;R))$$

For the proofs of the following three theorems one can see [1].

 {\bf Theorem 2.} The following equality
$$Exp(M_{\varphi}Exp(x))=Exp(M_{\varphi})Exp(x) $$ is valid.

 Consider
$\psi(y)=(\psi_1(y),\psi_2(y),...,\psi_{n'}(y))=M_{\psi}Exp{(y)}$,
where $ M_{\psi}(1,i)\in M_{n'',n'}(1,i;R)$ and

$$(\varphi \circ \psi)(y)=(\varphi_1(\psi(y)),
\varphi_2(\psi(y)),...,\varphi_{n}(\psi(y)))=M_{\varphi \circ
\psi}Exp{(y)}$$

The following result is about the matrix representation of the
composition $\varphi \circ \psi$.

{\bf Theorem 3.} The following equality $$M_{\varphi \circ
\psi}=M_{\varphi}Exp(M_{\psi})$$ is valid.

{\bf Theorem 4.} The following equality
$$Exp(M_{\varphi}Exp(M_{\psi}))=Exp(M_{\varphi})Exp(M_{\psi})$$ is valid.

{\bf The behavior of the matrix of a polynomial map with respect to the change of variables.}
Now let us investigate the change of matrix of a map $\varphi
:R^{n'}\rightarrow R^{n}$ with respect to, not compulsory
linear, change of variables (coordinate system) in $R^{n'}$ and
$R^{n}$. Assume that with respect to some coordinate systems in
$R^{n'}$ and $R^{n}$ one has
$$y=\varphi(x)=M_{\varphi,x,y}Exp(x) $$ and with respect to
some coordinate systems in $R^{n'}$ and $R^{n}$ one has
$$y'=\varphi(x')=M_{\varphi,x',y'}Exp(x') $$ Here we are using
the notation $M_{\varphi,x,y}$ to indicate that
$M_{\varphi,x,y}$ is the matrix of $\varphi$ with respect to the
first coordinate systems in $R^{n'}$ and $R^{n}$. If
$x'=T_{x,x'}Exp(x)$ and $y'=S_{y,y'}Exp(y)$, where $T_{x,x'}$,
$S_{y,y'}$ are the corresponding transformation matrices  of variables changes, then one has

$$M_{\varphi,x',y'}Exp(x')=y'=S_{y,y'}Exp(y)=S_{y,y'}Exp(M_{\varphi ,x,y}Exp(x))=S_{y,y'}Exp(M_{\varphi,x,y})Exp(x)=$$
 $$S_{y,y'}Exp(M_{\varphi ,x,y})Exp(T_{x',x}Exp(x'))=S_{y,y'}Exp(M_{\varphi ,x,y})Exp(T_{x',x})Exp(x')$$

so $$M_{\varphi,x',y'}=S_{y,y'}Exp(M_{\varphi ,x,y})Exp(T_{x',x})$$

If $\varphi(x)=A(1,k)\frac{x^{(k)}}{k!}$ (that is $\varphi$ is a
homogeneous map of degree $k$) and $x'=Tx$, $y'=Sy$ are
 linear transformations then $$A'(1,k)=SA(1,k)\frac{(T^{-1})^{(k)}}{k!}$$

 If in addition $n=n'$ and $S=T$ then one has
$$A'(1,k)=TA(1,k)\frac{(T^{-1})^{(k)}}{k!}$$
Therefore problem of equivalence of $"1\times k"$ size matrices
with respect to such action of the group $GL(n,R)$ can be
considered as an interesting one. One can see that rank of
$A(1,k)$ is one of the invariants with respect to this action.

{\bf Example 2.} If $n=k=2$ in ordinary notations for matrices
$A=\left(
\begin{array}{ccc}
  a_{11} & a_{12}& a_{13} \\
  a_{21} & a_{22} & a_{23} \\
  \end{array}
\right)$, \\ $T^{-1}=\left(
\begin{array}{cc}
  t_{11} & t_{12} \\
  t_{21} & t_{22} \\
\end{array}
\right)$ the above equality can be written as
$$\left(
\begin{array}{ccc}
  a'_{11} & a'_{12}& a'_{13} \\
  a'_{21} & a'_{22} & a'_{23} \\
  \end{array}
\right)=\left(
\begin{array}{cc}
  t_{11} & t_{12} \\
  t_{21} & t_{22} \\
\end{array}
\right)^{-1}\left(
\begin{array}{ccc}
  a_{11} & a_{12}& a_{13} \\
  a_{21} & a_{22} & a_{23} \\
  \end{array}
\right) \left(
\begin{array}{ccc}
  t^2_{11} & t_{11}t_{12}& t^2_{12} \\
  2t_{11}t_{21} & 2(t_{11}t_{22}+t_{12}t_{21}) & 2t_{12}t_{22} \\
  t^2_{21} & t_{21}t_{22} & t^2_{22}\\
\end{array}
\right)$$

{\bf Matrix representation for symmetric multi-linear maps.}
The above introduced product $\bigodot$ is convenient to represent
symmetric multi-linear maps by the use of matrices as well. If $V$
($V'$) is a $n$-dimensional (respect. $n'$-dimensional) vector
space with a given basis and $\textbf{A}:V^p\rightarrow V'$ is a
symmetric multi-linear map one can attach to $\textbf{A}$ a matrix
$A$ for which the equality
$$\textbf{A}(x^1,x^2,...,x^p)=A(\frac{x^1\bigodot x^2\bigodot ...\bigodot x^p}{p!})$$
is true at all $x^1,x^2,...,x^p$ from $V$. In other words to
$\textbf{A}$ is attached the same matrix $A$ which was attached to
the corresponding homogeneous polynomial map $\textbf{A}(x,x,...,x)$
in [1].

If $\textbf{B}:V^q\rightarrow V''$ is a symmetric multi-linear map
and $\textbf{C}:V'\times V''\rightarrow V'''$ is a given bilinear
map one can consider the "product"
$\textbf{A}\times_{\textbf{C}}\textbf{B}$ of $\textbf{A}$ and $\textbf{B}$
with respect to $\textbf{C}(x',x'')=C\frac{(x',x'')^{(2)}}{2!}$ as a
symmetric multi-linear map
$\textbf{A}\times_{\textbf{C}}\textbf{B}:V^{p+q}\rightarrow V'''$
defined by
$$\textbf{A}\times_{\textbf{C}}\textbf{B}(x^1,x^2,...,x^{p+q})=\frac{1}{(p+q)!}\sum_{\sigma\in
S_{p+q}}\textbf{C}(A\frac{x^{\sigma(1)}\bigodot
x^{\sigma(2)}\bigodot ...\bigodot
x^{\sigma(p)}}{p!},B\frac{x^{\sigma(p+1)}\bigodot
x^{\sigma(p+2)}\bigodot ...\bigodot x^{\sigma(p+q)}}{q!})$$, where
$S_{p+q}$ stands for the symmetric group.

{\bf Theorem 5.} The equality
$$\textbf{A}\times_{\textbf{C}}\textbf{B}(x^1,x^2,...,x^{p+q})=C(\overline{A}\bigodot
\underline{B})\frac{x^1\bigodot x^2\bigodot ...\bigodot
x^{p+q}}{(p+q)!}$$ is valid, where
$\overline{A}=\left(\begin{array}{c}
  A \\
  0 \\
\end{array}\right)\in M_{n'+n'',n}(1,p;R)$, $\underline{B}=\left(\begin{array}{c}
  0 \\
  B \\
\end{array}\right)\in M_{n'+n'',n}(1,q;R)$.

{\bf Proof.} The  map $\textbf{C}$ is a bilinear map and therefore
$\textbf{C}(x',0)=\textbf{C}(0,x'')=0$,
$\textbf{C}(x',x'')=C((x',0)\bigodot(0,x''))$. It implies that
$$\textbf{A}\times_{\textbf{C}}\textbf{B}(x^1,x^2,...,x^{p+q})=$$
$$\frac{1}{(p+q)!}\sum_{\sigma\in
S_{p+q}}C((A\frac{x^{\sigma(1)}\bigodot
x^{\sigma(2)}\bigodot ...\bigodot x^{\sigma(p)}}{p!},0)\bigodot
(0, B\frac{x^{\sigma(p+1)}\bigodot x^{\sigma(p+2)}\bigodot
...\bigodot x^{\sigma(p+q)}}{q!}))$$

To complete the proof it is enough to show that
$$\sum_{\sigma\in
S_{p+q}}(A\frac{x^{\sigma(1)}\bigodot x^{\sigma(2)}\bigodot
...\bigodot x^{\sigma(p)}}{p!},0)\bigodot (0,
B\frac{x^{\sigma(p+1)}\bigodot x^{\sigma(p+2)}\bigodot ...\bigodot
x^{\sigma(p+q)}}{q!})=$$ $$(\overline{A}\bigodot
\underline{B})(x^1\bigodot x^2\bigodot ...\bigodot x^{p+q})$$

For any $1\leq i\leq p$, $p< j\leq p+q$ the $(e_i+e_j,0)$
component, where $e_k$ means $e_k\in M_{0,p+q}(0,1;R)$ all
coordinates are zero except for $k$-th which is 1, of the left
side is equal to

$$\sum_{\sigma\in S_{p+q}}\sum_{|\alpha
|=p}A_{e_i,\alpha}(\frac{x^{\sigma(1)}\bigodot
x^{\sigma(2)}\bigodot ...\bigodot
x^{\sigma(p)}}{p!})_{\alpha,0}\sum_{|\beta|=q}
B_{e_j,\beta}(\frac{x^{\sigma(p+1)}\bigodot
x^{\sigma(p+2)}\bigodot ...\bigodot
x^{\sigma(p+q)}}{q!})_{\beta,0}=$$ $$ \sum_{|\alpha |=p,
|\beta|=q}A_{e_i,\alpha}B_{e_j,\beta}\sum_{\sigma\in
S_{p+q}}(\frac{x^{\sigma(1)}\bigodot x^{\sigma(2)}\bigodot
...\bigodot x^{\sigma(p)}}{p!})_{\alpha,0}
(\frac{x^{\sigma(p+1)}\bigodot x^{\sigma(p+2)}\bigodot ...\bigodot
x^{\sigma(p+q)}}{q!})_{\beta,0}$$

But by comparison of the corresponding terms it can be checked
that

$$\sum_{\sigma\in
S_{p+q}}(\frac{x^{\sigma(1)}\bigodot x^{\sigma(2)}\bigodot
...\bigodot
x^{\sigma(p)}}{p!})_{(\alpha,0)}(\frac{x^{\sigma(p+1)}\bigodot
x^{\sigma(p+2)}\bigodot ...\bigodot
x^{\sigma(p+q)}}{q!})_{(\beta,0)}=$$ $$\left(\begin{array}{c}
  \alpha+\beta \\
  \beta \\
\end{array}\right)(x^1\bigodot x^2\bigodot ...\bigodot
x^{p+q})_{(\alpha+\beta,0)}$$

It implies that the result is true.

{\bf Remark 1.}  If $\textbf{C}:V'\times V'\rightarrow V'''$ is itself a symmetric bilinear
map then one can attach to it a matrix $C$ for which $\textbf{C}(x',y')=C(\frac{x'\bigodot y'}{2!})$ is true. In this case if
$\textbf{A}:V^p\rightarrow V'$, $\textbf{B}:V^q\rightarrow V'$ are symmetric multi-linear maps then
$$\textbf{A}\times_{\textbf{C}}\textbf{B}(x^1,x^2,...,x^{p+q})=\frac{1}{2!}C(A\bigodot
B)\frac{x^1\bigodot x^2\bigodot ...\bigodot
x^{p+q}}{(p+q)!}$$

{\bf Matrix representation for antisymmetric multi-linear maps.}
Now let us consider antisymmetric multi-linear maps.
In future we use the following notations.
$$S_{[p]+[q]+[r]}=\{\sigma\in S_{p+q+r}: \sigma(1)<...<\sigma(p);\ \sigma(p+1)<...<\sigma(p+q); \ \sigma(p+q+1)<...<\sigma(p+q+r)\}$$
$$S_{[p]+[q+r]}=\{\sigma\in S_{p+q+r}: \sigma(1)<...<\sigma(p);\ \sigma(p+1)<...<\sigma(p+q+r)\}$$
$$S_{(p)+[q]+[r]}=\{\sigma\in S_{p+q+r}: \sigma(i)=i \ \mbox{for}\  i=\overline{1,p};\ \sigma(p+1)<...<\sigma(p+q); \ \sigma(p+q+1)<...<\sigma(p+q+r)\}$$

{\bf Proposition 4.} The map $P: S_{[p]+[q+r]}\times S_{(p)+[q]+[r]}\rightarrow S_{[p]+[q]+[r]}$, where $P(\sigma, \tau)=\sigma\tau$ is an one to one correspondence
between the sets $S_{[p]+[q+r]}\times S_{(p)+[q]+[r]}$ and $S_{[p]+[q]+[r]}$.

{\bf Proof.} It is easy to see that indeed $P$ maps $S_{[p]+[q+r]}\times S_{(p)+[q]+[r]}$ into $S_{[p]+[q]+[r]}$.
To show that $P$ is "on" for $\sigma_0\in S_{[p]+[q]+[r]}$ define $\sigma\in S_{[p]+[q+r]}$ as following: $$\sigma(i)=\sigma_0(i)\ \mbox{for}\ i=\overline{1,p}, [ \sigma_0(p+1), \sigma_0(p+2),..., \sigma_0(p+q+r)]= (\sigma(p+1), \sigma(p+2),..., \sigma(p+q+r))$$, where $[ \sigma_0(p+1), \sigma_0(p+2),..., \sigma_0(p+q+r)]$ means ordering the numbers $ \sigma_0(p+1), \sigma_0(p+2),..., \sigma_0(p+q+r)$ in ascending order. Now one can take for $\tau$ the element $\sigma^{-1}\sigma_0$ which is in $S_{(p)+[q]+[r]}$. To show that $P$ is "one to one" assume that $\sigma_0=\sigma\tau$, where $\sigma_0\in S_{[p]+[q]+[r]}$, $\sigma\in S_{[p]+[q+r]}$ and $\tau\in S_{(p)+[q]+[r]}$. In this case the $\sigma$ is unique as far as $$\sigma(i)=\sigma_0(i)\ \mbox{for}\ i=\overline{1,p}\ \mbox{and}\ [ \sigma_0(p+1), \sigma_0(p+2),..., \sigma_0(p+q+r)]= (\sigma(p+1), \sigma(p+2),..., \sigma(p+q+r))$$

This result is valid in the following form as well:

{\bf Proposition 4'.} The map $P': S_{[p+q]+[r]}\times S_{[p]+[q]+(r)}\rightarrow S_{[p]+[q]+[r]}$,
where $P'(\sigma, \tau)=\sigma\tau$ is an one to one correspondence
between the sets $S_{[p+q]+[r]}\times S_{[p]+[q]+(r)}$ and $S_{[p]+[q]+[r]}$.

 Let $n, n'$ be
any fixed nonnegative integers, $J_n(p)=\{\alpha
=(\alpha_1,\alpha_2,...,\alpha_p): 1\leq \alpha_1 < \alpha_2
<...<\alpha_p \leq n \}$, $J_n(0)=\{0\}$. If $\alpha\in J_n(p), \beta\in J_n(q)$ we say $\alpha < \beta$ if $p<q$, or $p=q$ and
$\alpha_p < \beta_p$, or $p=q$, $\alpha_p = \beta_p$ and $\alpha_{p-1} < \beta_{p-1}$ and etcetera.
We use the following notations as well. $$\textbf{M}(p,p';R)=\textbf{M}_{n,n'}(p,p';R)=\{A=A(p,p')=(A_{\alpha,\alpha'})_{\alpha \in J_n(p),
\alpha'\in J_{n'}(p')}: A_{\alpha,\alpha'}\in R \}$$

The ordinary size of $A(p,p')$ is $\left(\begin{array}{c}
  n \\
  p \\
\end{array}\right)\times \left(\begin{array}{c}
  n' \\
  p' \\
\end{array}\right) $.

{\bf Definition 2.} If $A=A(p,p')\in \textbf{M}(p,p';R)$ and $B=B(q,q')\in \textbf{M}(q,q';R)$ then
$A\wedge B=C\in \\  \textbf{M}(p+q, p'+q';R)$ such that for any $\alpha =(\alpha_1,
\alpha_2,...,\alpha_{p+q})\in J_n(p+q)$, $\alpha' =(\alpha'_1,
\alpha'_2,...,\alpha'_{p'+q'})\in J_{n'}(p'+q')$
$$C_{\alpha,\alpha'}=\sum\varepsilon(\sigma)\varepsilon(\sigma')
A_{(\alpha_{\sigma(1)},
...,\alpha_{\sigma(p)}),(\alpha'_{\sigma'(1)},
...,\alpha'_{\sigma'(p')})}B_{(\alpha_{\sigma(p+1)},
...,\alpha_{\sigma(p+q)}),(\alpha'_{\sigma'(p'+1)},
...,\alpha'_{\sigma'(p'+q')})}
$$, where the sum is taken over all $\sigma\in S_{[p]+[q]}$ and $\sigma'\in
S_{[p']+[q']}$.

{\bf Example 3.}  a) If $n=n'=2$, $p=p'=q=q'=1$ in ordinary
notations for matrices $A=\left(
\begin{array}{cc}
  a_{11} & a_{12} \\
  a_{21} & a_{22} \\
\end{array}
\right)$, $B=\left(
\begin{array}{cc}
  b_{11} & b_{12} \\
  b_{21} & b_{22} \\
\end{array}
\right)$ the product $A\wedge B$ can be given as
$$A\wedge B=(a_{11}b_{22}+ a_{22}b_{12}-a_{12}b_{21}-a_{21}b_{12})$$

b) If $n=n'=3$, $p=p'=q=q'=1$ in ordinary notations for matrices
$A=\left(
\begin{array}{ccc}
  a_{11} & a_{12}& a_{13} \\
  a_{21} & a_{22} & a_{23} \\
  a_{31} & a_{32}& a_{33} \\
  \end{array}
\right)$, $B=\left(
\begin{array}{ccc}
  b_{11} & b_{12}& b_{13} \\
  b_{21} & b_{22} & b_{23} \\
  b_{31} & b_{32}& b_{33} \\
  \end{array}
\right)$ the product $A\wedge B$, an ordinary $3\times 3$ size
matrix, can be given as
$$A\wedge B= \left(
\begin{array}{ccc}
 \left|\begin{array}{cc}
    a_{11} & a_{12} \\
    b_{21} & b_{22} \\
  \end{array}\right|+\left|\begin{array}{cc}
    b_{11} & b_{12} \\
    a_{21} & a_{22} \\
  \end{array}\right| &
  \left|\begin{array}{cc}
    a_{11} & a_{13} \\
    b_{21} & b_{23} \\
  \end{array}\right|+\left|\begin{array}{cc}
    b_{11} & b_{13} \\
    a_{21} & a_{23} \\
  \end{array}\right| &
  \left|\begin{array}{cc}
    a_{12} & a_{13} \\
    b_{22} & b_{23} \\
  \end{array}\right|+\left|\begin{array}{cc}
    b_{12} & b_{13} \\
    a_{22} & a_{23} \\
  \end{array}\right| \\
      . & . & .\\
     \left|\begin{array}{cc}
    a_{11} & a_{12} \\
    b_{31} & b_{32} \\
  \end{array}\right|+\left|\begin{array}{cc}
    b_{11} & b_{12} \\
    a_{31} & a_{32} \\
  \end{array}\right| &
  \left|\begin{array}{cc}
    a_{11} & a_{13} \\
    b_{31} & b_{33} \\
  \end{array}\right|+\left|\begin{array}{cc}
    b_{11} & b_{13} \\
    a_{31} & a_{33} \\
  \end{array}\right| &
  \left|\begin{array}{cc}
    a_{12} & a_{13} \\
    b_{32} & b_{33} \\
  \end{array}\right|+\left|\begin{array}{cc}
    b_{12} & b_{13} \\
    a_{32} & a_{33} \\
  \end{array}\right| \\
   . & . & .\\
  \left|\begin{array}{cc}
    a_{21} & a_{22} \\
    b_{31} & b_{32} \\
  \end{array}\right|+\left|\begin{array}{cc}
    b_{21} & b_{22} \\
    a_{31} & a_{32} \\
  \end{array}\right| &
  \left|\begin{array}{cc}
    a_{21} & a_{23} \\
    b_{31} & b_{33} \\
  \end{array}\right|+\left|\begin{array}{cc}
    b_{21} & b_{23} \\
    a_{31} & a_{33} \\
  \end{array}\right| &
  \left|\begin{array}{cc}
    a_{22} & a_{23} \\
    b_{32} & b_{33} \\
  \end{array}\right|+\left|\begin{array}{cc}
    b_{22} & b_{23} \\
    a_{32} & a_{33} \\
  \end{array}\right| \\
\end{array}\right)
  $$

{\bf Proposition 5.} For the above defined product the following
are true.

1. $A(p,p')\wedge B(q,q')=(-1)^{pq+p'q'}B(q,q')\wedge A(p,p')$.

2. $(A+B)\wedge C=A\wedge C+ B\wedge C$.

3. $(A\wedge B)\wedge C=A\wedge (B\wedge C)$

4. $ (\lambda A)\wedge B=\lambda (A\wedge B)$ for any
$\lambda\in R$

5. If $A$ and $B$ are square upper triangular matrices then $A\wedge B $
is also an upper triangular matrix.

{\bf Proof.} Let us prove the associative property: $$(A(p,p')\wedge
B(q,q'))\wedge C(r,r')=A(p,p')\wedge (B(q,q')\wedge C(r,r'))$$ To
show this equality one can consider $((\alpha_1,
\alpha_2,...,\alpha_{p+q+r}),(\alpha'_1,\alpha'_2,...,\alpha'_{p'+q'+r'}))$
component, where \\ $1\leq\alpha_1< ...<\alpha_{p+q+r}\leq n$,
$1\leq\alpha'_1<...<\alpha'_{p'+q'+r'}\leq n'$, of the left side
of it:
$$((A\wedge B)\wedge C)_{(\alpha_1, \alpha_2,...,\alpha_{p+q+r}),(\alpha'_1,\alpha'_2,...,\alpha'_{p'+q'+r'})}=$$
$$\sum_{\sigma\in S_{[p+q]+[r]},\sigma'\in S_{[p'+q']+[r']}}\varepsilon(\sigma)\varepsilon(\sigma')
(A\wedge B)_{(\alpha_{\sigma(1)},
...,\alpha_{\sigma(p+q)}),(\alpha'_{\sigma'(1)},
...,\alpha'_{\sigma'(p'+q')})}$$
$$C_{(\alpha_{\sigma(p+q+1)},
...,\alpha_{\sigma(p+q+r)}),(\alpha'_{\sigma'(p'+q'+1)},
...,\alpha'_{\sigma'(p'+q'+r')})}=$$
$$\sum_{\sigma\in S_{[p+q]+[r]},\sigma'\in S_{[p'+q']+[r']}}\varepsilon(\sigma)\varepsilon(\sigma')
\sum_{\tau\in S_{[p]+[q]+(r)},\tau'\in
S_{[p']+[q']+(r')}}\varepsilon(\tau)\varepsilon(\tau')
A_{(\alpha_{\sigma(\tau(1))},
...,\alpha_{\sigma(\tau(p))}),(\alpha'_{\sigma'(\tau'(1))},
...,\alpha'_{\sigma'(\tau'(p'))})}$$
$$B_{(\alpha_{\sigma(\tau(p+1))},
...,\alpha_{\sigma(\tau(p+q))}),(\alpha'_{\sigma'(\tau'(p'+1))},
...,\alpha'_{\sigma'(\tau'(p'+q'))})}C_{(\alpha_{\sigma(p+q+1)},
...,\alpha_{\sigma(p+q+r)}),(\alpha'_{\sigma'(p'+q'+1)},
...,\alpha'_{\sigma'(p'+q'+r')})}$$ Due to Proposition 4' the
right side of the last equality can be written in the form
$$\sum_{\sigma\in S_{[p]+[q]+[r]},\sigma'\in S_{[p']+[q']+[r']}}\varepsilon(\sigma)\varepsilon(\sigma')
A_{(\alpha_{\sigma(1)},
...,\alpha_{\sigma(p)}),(\alpha'_{\sigma'(1)},
...,\alpha'_{\sigma'(p')})}B_{(\alpha_{\sigma(p+1)},
...,\alpha_{\sigma (p+q))}),(\alpha'_{\sigma'(p'+1)},
...,\alpha'_{\sigma'(p'+q')})}$$
$$C_{(\alpha_{\sigma(p+q+1)},
...,\alpha_{\sigma(p+q+r)}),(\alpha'_{\sigma'(p'+q'+1)},
...,\alpha'_{\sigma'(p'+q'+r')})}$$ In similar way due to
Proposition 4 for the right side of the initial equality to be proved one will have the same
expression.

{\bf Proposition 6.} For any $1\leq\alpha^1_1 <...<\alpha^1_{p_1} < \alpha^2_{p_1+1}<...<\alpha^2_{p_1+p_2}<...< \alpha^k_{p_1+p_2+...+p_k}\leq n $ and
$1\leq\alpha'^1_1 <...<\alpha'^1_{p'_1} < \alpha'^2_{p'_1+1}<...<\alpha'^2_{p'_1+p'_2}<...< \alpha'^k_{p'_1+p'_2+...+p'_k}\leq n'$
the following equality is true
$$(A^1(p_1,p'_1)\wedge A^2(p_2,p'_2)\wedge ...\wedge A^k(p_k,p'_k))_{(\alpha^1_1,...,\alpha^1_{p_1},\alpha^2_{p_1+1},...,\alpha^2_{p_1+p_2},..., \alpha^k_{p_1+p_2+...+p_k}),(\alpha'^1_1,...,\alpha^1_{p'_1},..., \alpha^k_{p'_1+p'_2+...+p'_k})}=$$ $$
\sum \varepsilon(\sigma)\varepsilon(\sigma')A^1_{(\alpha^1_{\sigma(1)},...,\alpha^1_{\sigma(p_1)}),(\alpha'^1_{\sigma'(1)},...,\alpha'^1_{\sigma'(p'_1)})}
A^2_{(\alpha^2_{\sigma(p_1+1)},...,\alpha^2_{\sigma(p_1+p_2)}),(\alpha'^2_{\sigma'(p'_1+1)},...,\alpha'^2_{\sigma'(p'_1+p'^2)})}...$$ $$
A^k_{(\alpha^k_{\sigma(p_1+...+p_{k-1}+1)},...,\alpha^k_{\sigma(p_1+...+p_k)}),(\alpha'^k_{\sigma'(p'_1+...+p'_{k-1}+1)},...,\alpha'^k_{\sigma'(p'_1+...+p'_k)})}$$
, where the sum is taken over all $\sigma\in S_{[p_1]+[p_2]+...+[p_k]}$ and $\sigma'\in S_{[p'_1]+[p'_2]+...+[p'_k]}$.

This proposition can be proved by induction on $k$ taking into account the associative property of the "product" $\wedge$.

{\bf Theorem 6.} If $A$ is a $"1\times 1"$ size matrix, $1\leq\alpha_1 <...<\alpha_{k} \leq n $ and
$1\leq\alpha'_1 <...<\alpha'_k \leq n'$ then

$1.\ \ (A\wedge A\wedge ...\wedge A)_{(\alpha_1,...,\alpha_{k}),(\alpha'_1,...,\alpha'_{k})}=A^{\wedge k}_{(\alpha_1,...,\alpha_{k}),(\alpha'_1,...,\alpha'_{k})}=$\\
$
k!\sum_{\sigma\in S_k} \varepsilon(\sigma)A_{\alpha_{\sigma(1)},\alpha'_1}A_{\alpha_{\sigma(2)},\alpha'_2}...
A_{\alpha_{\sigma(k)},\alpha'_k}=k!\sum_{\sigma\in S_k} \varepsilon(\sigma)A_{\alpha_1,\alpha'_{\sigma(1)}}A_{\alpha_2,\alpha'_{\sigma(2)}}...
A_{\alpha_k,\alpha'_{\sigma(k)}} $

$2.\ \ Ax^1\wedge Ax^2\wedge ...\wedge Ax^k=\frac{A^{\wedge
k}}{k!}(x^1\wedge x^2\wedge...\wedge x^k)$

3. \ \ If $Av_i=\lambda_iv_i$ for $i=\overline{1,k}$ and $v_1,
v_2,..., v_k$ are linear independent then $v_1\wedge
v_2\wedge...\wedge v_k$ is an eigenvector for $\frac{A^{\wedge
k}}{k!}$ with the eigenvalue $\lambda_1\lambda_2...\lambda_k$.

Let $A$ be a $"1\times 1"$ size square matrix.

4.  If $\{\lambda_1,\lambda_2,...,\lambda_n\}$ is all its
eigenvalues, where each eigenvalue occurs as many times as its
multiplicity, then for $\frac{A^{\wedge k}}{k!}$ the system
$\{\lambda_{i_1}\lambda_{i_2}...\lambda_{i_k}\}_{1\leq
i_1<....<i_k\leq n}$ represents all its eigenvalues including
their multiplicities as eigenvalues of  $\frac{A^{\wedge k}}{k!}$.

5. If $rk(A)=l$ then $rk(A^{\wedge k})=\left(\begin{array}{c}
  l \\
  k \\
\end{array}\right)$

6. $\det (\frac{A^{\wedge k}}{k!})=\det(A)^{\left(\begin{array}{c}
  n-1 \\
  k-1 \\
\end{array}\right)}$

{\bf Proof.} 1. For each $\sigma'\in S_k$ one has $\sum_{\sigma\in S_k} \varepsilon(\sigma) \varepsilon(\sigma') A_{\alpha_{\sigma(1)},\alpha'_{\sigma'(1)}}A_{\alpha_{\sigma(2)},\alpha'_{\sigma'(2)}}...
A_{\alpha_{\sigma(k)},\alpha'_{\sigma'(k)}}=$ $$\sum_{\sigma\in S_k} \varepsilon(\sigma) \varepsilon(\sigma') A_{\alpha_{\sigma(\sigma'^{-1}(1))},\alpha'_1}A_{\alpha_{\sigma(\sigma'^{-1}(2))},\alpha'_2}...
A_{\alpha_{\sigma(\sigma'^{-1}(k))},\alpha'_k}=\sum_{\sigma\in S_k} \varepsilon(\sigma)A_{\alpha_{\sigma(1)},\alpha'_1}A_{\alpha_{\sigma(2)},\alpha'_2}...
A_{\alpha_{\sigma(k)},\alpha'_k}$$ and therefore $$A^{\wedge k}_{(\alpha_1,...,\alpha_{k}),(\alpha'_1,...,\alpha'_{k})}=
\sum_{\sigma'\in S_k}\sum_{\sigma\in S_k} \varepsilon(\sigma)\varepsilon(\sigma')A_{\alpha_{\sigma(1)},\alpha'_1}A_{\alpha_{\sigma(2)},\alpha'_2}...
A_{\alpha_{\sigma(k)},\alpha'_k}=$$ $$\sum_{\sigma'\in S_k}\sum_{\sigma\in S_k} \varepsilon(\sigma)A_{\alpha_{\sigma(1)},\alpha'_1}A_{\alpha_{\sigma(2)},\alpha'_2}...
A_{\alpha_{\sigma(k)},\alpha'_k}=k!\sum_{\sigma\in S_k} \varepsilon(\sigma)A_{\alpha_{\sigma(1)},\alpha'_1}A_{\alpha_{\sigma(2)},\alpha'_2}...
A_{\alpha_{\sigma(k)},\alpha'_k}$$

2. For any $1\leq\alpha_1 <...<\alpha_{k} \leq n $  due to Proposition 6 one has $$(Ax^1\wedge Ax^2\wedge ...\wedge Ax^k)_{(\alpha_1,...,\alpha_{k}),0}=
\sum_{\sigma\in S_k} \varepsilon(\sigma)(Ax^1)_{\alpha_{\sigma(1)},0}(Ax^2)_{\alpha_{\sigma(2)},0}...
(Ax^k)_{\alpha_{\sigma(k)},0}=$$ $$\sum_{\sigma\in S_k} \varepsilon(\sigma)\sum_{s_1=\overline{1,n'}}A_{\alpha_{\sigma(1)},s_1}x^1_{s_1}...\sum_{s_k=
\overline{1,n'}}A_{\alpha_{\sigma(k)},s_k}x^k_{s_k}=
\sum_{s_1=\overline{1,n'},...,s_k=
\overline{1,n'}}x^1_{s_1}...x^k_{s_k}\sum_{\sigma\in S_k} \varepsilon(\sigma)A_{\alpha_{\sigma(1)},s_1}...A_{\alpha_{\sigma(k)},s_k}$$
If $s_1, s_2,...,s_k$ are not different then the last sum is equal to zero.
If $s_1, s_2,...,s_k$ are fixed, they are different numbers and $[s_1, s_2,...,s_k]=(s_{\tau_0(1)}, s_{\tau_0(2)},...,s_{\tau_0(k)}$
then for any $\tau\in S_k$ one has $[s_{\tau(1)}, s_{\tau(2)},...,s_{\tau(k)}]=(s_{\tau_0(1)}, s_{\tau_0(2)},...,s_{\tau_0(k)})$.
Moreover $$\sum_{\sigma\in S_k} \varepsilon(\sigma)A_{\alpha_{\sigma(1)},s_1}...A_{\alpha_{\sigma(k)},s_k}=
\sum_{\sigma\in S_k} \varepsilon(\sigma)A_{\alpha_{\sigma(\tau_0(1))},s_{\tau_0(1)}}...A_{\alpha_{\sigma(\tau_0(k))},s_{\tau_0(k)}}=$$
$$\varepsilon(\tau_0)\sum_{\sigma\in S_k} \varepsilon(\sigma)A_{\alpha_{\sigma(1)},s_{\tau_0(1)}}...A_{\alpha_{\sigma(k)},s_{\tau_0(k)}}=
\frac{\varepsilon(\tau_0)}{k!}A^{\wedge
k}_{(\alpha_1,...,\alpha_{k}),(s_{\tau_0(1)},...,s_{\tau_0(k)})}$$
So $$(Ax^1\wedge Ax^2\wedge ...\wedge
Ax^k)_{(\alpha_1,...,\alpha_{k}),0}=$$ $$\sum_{1\leq
s_1<...<s_k\leq n'} \frac{1}{k!}A^{\wedge
k}_{(\alpha_1,...,\alpha_{k}),(s_1,...,s_k)}\sum_{\sigma\in S_k}
\varepsilon(\sigma)x^{\sigma(1)}_{s_1}...x^{\sigma(k)}_{s_k}=(\frac{A^{\wedge
k}}{k!}(x^1\wedge x^2\wedge...\wedge
x^k))_{(\alpha_1,...,\alpha_{k}),0}$$

3. It is an easy consequence of the second statement.

4. From the second statement of this theorem one can conclude that
$$\frac{(T^{-1}AT)^{\wedge k}}{k!}=\frac{(T^{-1})^{\wedge k}}{k!}\frac{A^{\wedge
k}}{k!}\frac{T^{\wedge k}}{k!}$$, where $T$ is any $"1\times 1"$
size invertible matrix. Therefore one can choose $T$ to make
$T^{-1}AT$ in Jordan form. In this case due to Proposition 5
$\frac{(T^{-1}AT)^{\wedge k}}{k!}$ will be an upper triangular
matrix with $\{\lambda_{i_1}\lambda_{i_2}...\lambda_{i_k}\}_{1\leq
i_1<....<i_k\leq n}$ on the main diagonal. Therefore the fourth
statement of the theorem is also true.

5. It can be easily derived from the fourth statement.

6. It can be easily derived from the fourth statement.

If $V$ ($V'$) is a $n$-dimensional (respect. $n'$-dimensional)
vector space with a given basis and $\textbf{A}:V^{p}\rightarrow V'$
is an antisymmetric multi-linear map one can attach to $\textbf{A}$ a
matrix $A$ for which the equality
$$\textbf{A}(x^1,x^2,...,x^p)=A(x^1\wedge x^2\wedge ...\wedge x^p)$$
is true at all $x^1,x^2,...,x^p$ from $V$.

{\bf Remark 2.} If $V=V'$ and one changes the basis of $V$ he gets
for the new matrix $A'$ of the antisymmetric multi-linear map
$\textbf{A}$ an expression like $$A'=T^{-1}A \frac{T^{\wedge
p}}{p!}$$, where $T\in GL(n,R)$. Therefore investigation of such
actions of $GL(n,R)$ also is an interesting problem. In the case
$n=3$, $p=2$ in ordinary notations for matrices this action can be
expressed as
$$A'=\left(
\begin{array}{ccc}
  t_{11} & t_{12} & t_{13}\\
  t_{21} & t_{22} & t_{23}\\
  t_{31} & t_{32} & t_{33}\\
\end{array}
\right)^{-1}\left(
\begin{array}{ccc}
  a_{11} & a_{12} & a_{13}\\
  a_{21} & a_{22} & a_{23}\\
  a_{31} & a_{32} & a_{33}\\
\end{array}
\right)\left(
\begin{array}{ccc}
 \left|\begin{array}{cc}
    t_{11} & t_{12} \\
    t_{21} & t_{22} \\
  \end{array}\right| &
  \left|\begin{array}{cc}
    t_{11} & t_{13} \\
    t_{21} & t_{23} \\
  \end{array}\right|&
  \left|\begin{array}{cc}
    t_{12} & t_{13} \\
    t_{22} & t_{23} \\
  \end{array}\right| \\
      . & . & .\\
    \left|\begin{array}{cc}
    t_{11} & t_{12} \\
    t_{31} & t_{32} \\
  \end{array}\right|&
  \left|\begin{array}{cc}
    t_{11} & t_{13} \\
    t_{31} & t_{33} \\
  \end{array}\right|&
  \left|\begin{array}{cc}
    t_{12} & t_{13} \\
    t_{32} & t_{33} \\
  \end{array}\right| \\
   . & . & .\\
  \left|\begin{array}{cc}
    t_{21} & t_{22} \\
    t_{31} & t_{32} \\
  \end{array}\right|&
  \left|\begin{array}{cc}
    t_{21} & t_{23} \\
    t_{31} & t_{33} \\
  \end{array}\right| &
  \left|\begin{array}{cc}
    t_{22} & t_{23} \\
    t_{32} & t_{33} \\
    \end{array}\right| \\
\end{array}\right)  $$

Due to [3] if $\textbf{B}:V^q\rightarrow V''$ is an antisymmetric multi-linear
map and $\textbf{C}:V'\times V''\rightarrow V'''$ is a given
bilinear map one can consider the "product"
$\textbf{A}\wedge_{\textbf{C}}\textbf{B}$ of $\textbf{A}$ and $\textbf{B}$
with respect to $\textbf{C}(x',x'')=C\frac{(x',x'')^{(2)}}{2!}$ as an
antisymmetric multi-linear map
$\textbf{A}\wedge_{\textbf{C}}\textbf{B}:V^{p+q}\rightarrow V'''$
defined by
$$\textbf{A}\wedge_{\textbf{C}}\textbf{B}(x^1,x^2,...,x^{p+q})=\sum_{\sigma\in S_{[p]+[q]}}\varepsilon (\sigma)\textbf{C}(\textbf{A}(x^{\sigma(1)},
x^{\sigma(2)}, ..., x^{\sigma(p)}),\textbf{B}(x^{\sigma(p+1)},
x^{\sigma(p+2)}, ... , x^{\sigma(p+q)}))$$

{\bf Theorem 7.} The equality
$$\textbf{A}\wedge_{\textbf{C}}\textbf{B}(x^1,x^2,...,x^{p+q})=C(\overline{A}\wedge
\underline{B})(x^1\wedge x^2\wedge ...\wedge x^{p+q})$$ is valid,
where $\overline{A}=\left(\begin{array}{c}
  A \\
  0 \\
\end{array}\right)\in \textbf{M}_{n'+n'',n}(1,p;R)$, $\underline{B}=\left(\begin{array}{c}
  0 \\
  B \\
\end{array}\right)\in \textbf{M}_{n'+n'',n}(1,q;R)$.

{\bf Proof. } As far as
$\textbf{C}(x',x'')=C((x',0)\bigodot(0,x''))$ one has
$$\textbf{A}\times_{\textbf{C}}\textbf{B}(x^1,x^2,...,x^{p+q})=$$
$$\sum_{\sigma\in
S_{[p]+[q]}}\varepsilon (\sigma)C((A(x^{\sigma(1)}\wedge
x^{\sigma(2)}\wedge ...\wedge x^{\sigma(p)}),0)\bigodot
(0,B(x^{\sigma(p+1)}\wedge x^{\sigma(p+2)}\wedge ...\wedge
x^{\sigma(p+q)}))$$

To complete the proof it is enough to show that
$$\sum_{\sigma\in
S_{[p]+[q]}}\varepsilon (\sigma)(A(x^{\sigma(1)}\wedge
x^{\sigma(2)}\wedge ...\wedge x^{\sigma(p)}),0)\bigodot (0,
B(x^{\sigma(p+1)}\wedge x^{\sigma(p+2)}\wedge ...\wedge
x^{\sigma(p+q)}))=$$ $$(\overline{A}\wedge \underline{B})(x^1\wedge
x^2\wedge ...\wedge x^{p+q})$$

For any $1\leq i\leq p$, $p< j\leq p+q$ the $((i,j),0)$
component  of the left side is equal to

$$\sum_{\sigma\in S_{[p]+[q]}}\varepsilon
(\sigma)\sum_{\alpha'}A_{(i),\alpha'}(x^{\sigma(1)}\wedge
x^{\sigma(2)}\wedge ...\wedge
x^{\sigma(p)})_{\alpha',0}\sum_{\alpha''}
B_{(j),\alpha''}(x^{\sigma(p+1)}\wedge x^{\sigma(p+2)}\wedge
...\wedge x^{\sigma(p+q)})_{\alpha'',0}=$$
 $$ \sum_{\alpha',\alpha''}A_{(i),\alpha'}B_{(j),\alpha''}\sum_{\sigma\in
S_{[p]+[q]}}\varepsilon (\sigma)(x^{\sigma(1)}\wedge
x^{\sigma(2)}\wedge ...\wedge x^{\sigma(p)})_{\alpha',0}
(x^{\sigma(p+1)}\wedge x^{\sigma(p+2)}\wedge ...\wedge
x^{\sigma(p+q)})_{\alpha'',0}$$, where
$\alpha'=(\alpha_1,...,\alpha_p)$, $1\leq\alpha_1<...<\alpha_p\leq
n$, $\alpha''=(\alpha_{p+1},...,\alpha_{p+q})$,
$1\leq\alpha_{p+1}<...<\alpha_{p+q}\leq n$

Due to the equality

$$\det\left(\begin{array}{cccc}x^1_{1}&x^2_{1}&...&x^{p+q}_{1}\\
x^1_{2}&x^2_{2}&...&x^{p+q}_{2}\\
.&.&...&.\\.&.&...&.\\.&.&...&.\\
x^1_{p+q}&x^2_{p+q}&...&x^{p+q}_{p+q}\\
\end{array}\right)=(x^1\wedge
x^2\wedge ...\wedge x^{p+q})_{(1,2,...,p+q),0}= $$
$$\sum_{\sigma\in S_{[p]+[q]}}\varepsilon (\sigma)(x^{\sigma
(1)}\wedge x^{\sigma (2)}\wedge ...\wedge x^{\sigma
(p)})_{(1,2,...,p),0}(x^{\sigma (p+1)}\wedge x^{\sigma
(p+2)}\wedge ...\wedge x^{\sigma (p+q)})_{(p+1,p+2,...,p+q),0}=$$
$$\sum_{\sigma\in S_{[p]+[q]}}\varepsilon (\sigma)(x^1\wedge
x^2\wedge ...\wedge x^p)_{(\sigma (1),\sigma (2),...,\sigma
(p)),0}(x^{p+1}\wedge x^{p+2}\wedge ...\wedge
x^{p+q})_{(\sigma (p+1),\sigma (p+2),...,\sigma
(p+q)),0}$$ the sum $$\sum_{\sigma\in S_{[p]+[q]}}\varepsilon
(\sigma)(x^{\sigma(1)}\wedge x^{\sigma(2)}\wedge ...\wedge
x^{\sigma(p)})_{\alpha',0} (x^{\sigma(p+1)}\wedge
x^{\sigma(p+2)}\wedge ...\wedge x^{\sigma(p+q)})_{\alpha'',0}$$ is
zero whenever
$\{\alpha_1,...,\alpha_p\}\bigcap\{\alpha_{p+1},...,\alpha_{p+q}\}\neq
\emptyset $ and otherwise $$\sum_{\sigma\in S_{[p]+[q]}}\varepsilon
(\sigma)(x^{\sigma(1)}\wedge x^{\sigma(2)}\wedge ...\wedge
x^{\sigma(p)})_{\alpha',0} (x^{\sigma(p+1)}\wedge
x^{\sigma(p+2)}\wedge ...\wedge x^{\sigma(p+q)})_{\alpha'',0}=$$
$$\varepsilon(\tau)(x^1\wedge x^2\wedge ...\wedge
x^{p+q})_{(\alpha_{\tau(1)},...,\alpha_{\tau(p)},\alpha_{\tau(p+1)},...,\alpha_{\tau(p+q)}),0}$$,
where
$[\alpha_1,...,\alpha_p,\alpha_{p+1},...,\alpha_{p+q}]=(\alpha_{\tau(1)},...,\alpha_{\tau(p)},\alpha_{\tau(p+1)},...,\alpha_{\tau(p+q)})$

Therefore if
$1\leq\alpha_1<...<\alpha_p<\alpha_{p+1}<...<\alpha_{p+q}\leq n$
and one collects coefficients at \\ $(x^1\wedge x^2\wedge
...\wedge x^{p+q})_{\alpha,0}$, where $\alpha
=(\alpha_1,...,\alpha_p,\alpha_{p+1},...,\alpha_{p+q})$, in $$
\sum_{\alpha',\alpha''}A_{(i),\alpha'}B_{(j),\alpha''}\sum_{\sigma\in
S_{[p]+[q]}}\varepsilon (\sigma)(x^{\sigma(1)}\wedge
x^{\sigma(2)}\wedge ...\wedge x^{\sigma(p)})_{\alpha',0}
(x^{\sigma(p+1)}\wedge x^{\sigma(p+2)}\wedge ...\wedge
x^{\sigma(p+q)})_{\alpha'',0}$$ it will be equal to
$$\sum_{\sigma\in S_{[p]+[q]}}\varepsilon
(\sigma)A_{(i),(\alpha_{\sigma(1)},...,\alpha_{\sigma(p)})}B_{(j),(\alpha_{\sigma(p+1)},...,\alpha_{\sigma(p+q)})}=(\overline{A}\wedge
\underline{B})_{(i,j),\alpha}$$ It means that the matrix of
$\textbf{A}\wedge_{\textbf{C}}\textbf{B}$ is equal to $C(\overline{A}\wedge
\underline{B})$.

{\bf Remark 3.}  If $\textbf{C}:V'\times V'\rightarrow V'''$ is itself an antisymmetric bilinear
map then one can attach to it a matrix $C$ for which $\textbf{C}(x',y')=C(x'\wedge y')$ is true. In this case if
$\textbf{A}:V^p\rightarrow V'$, $\textbf{B}:V^q\rightarrow V'$ are antisymmetric multi-linear maps then
$$\textbf{A}\wedge_{\textbf{C}}\textbf{B}(x^1,x^2,...,x^{p+q})=C(A\wedge
B)(x^1\wedge x^2\wedge ...\wedge
x^{p+q})$$

If $\rho \geq 1$ is any real number and $F$ is the field of real
or complex numbers one can consider the following  norm in
$\textbf{M}_{n,n'}(p,p';F)$

{\bf Definition 3.} $$\|A\|_{\rho}=\|A\|=
(\sum_{\alpha\in J_n(p),\alpha'\in J_{n'}(p')}\frac{|A_{\alpha ,\alpha'}|^{\rho}}{(p!p'!)^{\rho-1}})^{\frac{1}{\rho}}$$

{\bf Theorem 8.} 1. If $A,B\in \textbf{M}(p,p';F)$ and $\lambda \in
F$ then

a)$\|A\|=0$ if and only if $A=0$,

b)$\|\lambda A\| =|\lambda |\|A\|$,

c) $\|A+B\|\leq \|A\|+ \|B\|$.

2. If $A\in \textbf{M}(p,p';F)$and $B\in \textbf{M}(q,q';F)$ then $\|A\wedge B\|\leq \|A\|\|B\|$

{\bf Proof.} Let $\varrho$ stand for the number for which $\frac{1}{\rho}+\frac{1}{\varrho}=1$
and let us prove only the second statement of the theorem. For any $1\leq\alpha_1< ...<\alpha_{p+q}\leq n$,
$1\leq\alpha'_1<...<\alpha'_{p'+q'}\leq n'$
 due to the H\"{o}lder inequality one has
 $$|(A\wedge B)_{\alpha,\alpha'}|=|\sum_{\sigma,\sigma'}(\varepsilon
(\sigma\sigma'))(A_{(\alpha_{\sigma(1)},
...,\alpha_{\sigma(p)}),(\alpha'_{\sigma'(1)},
...,\alpha'_{\sigma'(p')})}B_{(\alpha_{\sigma(p+1)},
...,\alpha_{\sigma(p+q)}),(\alpha'_{\sigma'(p'+1)},
...,\alpha'_{\sigma'(p'+q')})})|\leq $$
$$(\sum_{\sigma,\sigma'}1)^{\frac{1}{\varrho}}(\sum_{\sigma,\sigma'}|A_{(\alpha_{\sigma(1)},
...,\alpha_{\sigma(p)}),(\alpha'_{\sigma'(1)},
...,\alpha'_{\sigma'(p')})}B_{(\alpha_{\sigma(p+1)},
...,\alpha_{\sigma(p+q)}),(\alpha'_{\sigma'(p'+1)},
...,\alpha'_{\sigma'(p'+q')})}|^{\rho})^{\frac{1}{\rho}}=$$
$$(\left(\begin{array}{c}
  p+q \\
  p\\
\end{array}\right)\left(\begin{array}{c}
  p'+q' \\
  p' \\
\end{array}\right))^{\frac{1}{\varrho}}(\sum_{\sigma,\sigma'}
|A_{(\alpha_{\sigma(1)},
...,\alpha_{\sigma(p)}),(\alpha'_{\sigma'(1)},
...,\alpha'_{\sigma'(p')})}B_{(\alpha_{\sigma(p+1)},
...,\alpha_{\sigma(p+q)}),(\alpha'_{\sigma'(p'+1)},
...,\alpha'_{\sigma'(p'+q')})}|^{\rho})^{\frac{1}{\rho}}$$ as far as the number of elements of  $S_{[p]+[q]}$ is $\left(\begin{array}{c}
  p+q \\
  p\\
\end{array}\right)$, where the sums are taken over all \\ $\sigma\in S_{[p]+[q]},\sigma'\in S_{[p']+[q']}$. Therefore $$\|A\wedge B\|^{\rho}=\sum_{\alpha,\alpha'}\frac{|A_{\alpha ,\alpha'}|^{\rho}}{((p+q)!(p'+q')!)^{\rho-1}}\leq $$
$$(\frac{1}{p!p'!q!q'!})^{\rho -1}\sum_{\alpha,\alpha'}\sum_{\sigma,\sigma'}|A_{(\alpha_{\sigma(1)},
...,\alpha_{\sigma(p)}),(\alpha'_{\sigma'(1)},
...,\alpha'_{\sigma'(p')})}B_{(\alpha_{\sigma(p+1)},
...,\alpha_{\sigma(p+q)}),(\alpha'_{\sigma'(p'+1)},
...,\alpha'_{\sigma'(p'+q')})}|^{\rho}\leq \|A\|^{\rho}\|B\|^{\rho}$$ as far as

$$ \sum_{\alpha,\alpha'}\sum_{\sigma,\sigma'}|A_{(\alpha_{\sigma(1)},
...,\alpha_{\sigma(p)}),(\alpha'_{\sigma'(1)},
...,\alpha'_{\sigma'(p')})}B_{(\alpha_{\sigma(p+1)},
...,\alpha_{\sigma(p+q)}),(\alpha'_{\sigma'(p'+1)},
...,\alpha'_{\sigma'(p'+q')})}|^{\rho}\leq $$
 $$\sum_{\beta\in J_n{(p)},\beta'\in J_{n'}(p')}|A_{\beta,\beta'}|^{\rho}\sum_{\gamma\in J_n{(q)},\gamma'\in J_{n'}(q')}|B_{\gamma,\gamma'}|^{\rho}$$

\begin{center}{References}\end{center}

[1] Ural Bekbaev. \emph{A matrix representation of composition of
polynomial maps}. \\ arXiv0901.3179v3 [math. AC] 22 Sep.2009

[2]  Ural Bekbaev. \emph{A radius of absolute convergence for
power series in many variables}.\\ arXiv1001.0622v1 [math.CV] 5Jan.2010

[3] H.Cartan. Calcul diff\'{e}rentiel. Formes diff\'{e}rentielles.
Hermann. Paris, 1967.

\end{document}